\documentclass[10pt]{article}

\usepackage[latin1]{inputenc}
\usepackage{amsmath}
\usepackage{amsfonts,amssymb,latexsym,multicol}
\usepackage[francais]{babel}
\usepackage{epsfig,epsf}
\newcounter{fig}
\def\figcaption #1
  {\addtocounter{fig}{1}
   \begin{center}
   {\sc Fig.} \arabic{fig}. {\it #1} 
    \end{center}
    \addtocontents{lof}{\protect \contentsline {figure}{\protect \numberline {\thefig}{\protect \ignorespaces \protect \pit #1 }}{\thepage}}
   }
\input{cyracc.def}

\newtheorem{theo}{Th\'eor\`eme}

\newtheorem{prop}{Proposition}

\newcommand{\ioe}{\leqslant}
\newcommand{\soe}{\geqslant}
\newcommand{\vers}{\rightarrow}

\newcommand{\rad}{{\rm r}}

\newcommand{\fin}{\hfill$\Box$}
\newcommand{\dem}{\noindent {\bf D\'emonstration\ }}

\title{Sur l'infimum des parties r\'eelles des z\'eros des sommes partielles de la fonction z\^eta de Riemann}
\author{Michel Balazard et Oswaldo Vel\'asquez Casta\~n\'on}

\begin{document}
\maketitle

\begin{center}
  Let $ \varphi_n = \inf \bigl\{ \Re s \mid \sum_{m=1}^n m^{-s} = 0\}$. 
We show that $ \lim_{n \to +\infty} \varphi_n/n = - \log 2 $. 
\end{center}

\section{\'Enonc\'e du th\'eor\`eme principal et principes de la d\'emonstration}

Soit
$$ \zeta_n(s) = \sum_{m=1}^n m^{-s} \quad ( s=\sigma +i\tau \, ; \, \sigma ,\tau \in \mathbb{R})
$$
la somme partielle d'ordre $n \in  \mathbb{N}^*$  de la s\'erie de Dirichlet de la fonction z\^eta de Riemann. Posons
$$  \varphi_n = \inf \bigl\{ \sigma \mid \zeta_n(s) = 0\} \quad (n \in  \mathbb{N}^*).$$

\begin{theo} \label{varphigral} On a
$$ \lim_{n \to +\infty} \frac{\varphi_n}{n} = - \log 2. $$
\end{theo}

Pour $0 < k <n$, posons 
$$
\zeta_{n,k}(s)=n^{-s}-\sum_{n-k\ioe j<n}j^{-s}+\sum_{1\ioe j<n-k}j^{-s} \quad (s \in \mathbb{C}),
$$
et
$$  \rho_{n,k} = \inf \bigl\{ \sigma \mid \zeta_{n,k}(\sigma) = 0\} \quad (n \in  \mathbb{N}^*)
$$
(on ne consid\`ere ici que les z\'eros \emph{r\'eels}). Observons d\`es maintenant que $k\mapsto \rho_{n,k}$ est d\'ecroissante (car $k\mapsto \zeta_{n,k}(\sigma)$ est d\'ecroissante pour tout $\sigma \in \mathbb{R}$, et $\lim_{\sigma \to -\infty} \zeta_{n,k}(\sigma)=+\infty$).

Le th\'eor\`eme r\'esulte de la proposition suivante.

\begin{prop}
  (i) Pour tout $n \in  \mathbb{N}^*$, on a  
$\varphi_n \soe \rho_{n,n-1}$ ; (ii) Pour tout $k \in  \mathbb{N}^*$, il existe  $n_1(k) \in  \mathbb{N}^*$ tel que $\varphi_n \ioe \rho_{n,k}$ pour tout $n \soe n_1(k)$ ; (iii) 
$\rho_{n,k}/n \vers -\log 2 \quad (n>k \vers \infty).$
 \end{prop}

Le point \textit{(i)} est un exercice facile, signal\'e par exemple dans \cite[Theorem 3.1]{borwein}. Les auteurs de \cite{borwein} d\'emontrent en outre (\cite[p. 25]{borwein}) que
$$\frac{\rho_{n,n-1}}{n} \to -\log 2 \quad (n \to \infty).$$
Le  point \textit{(iii)} se d\'emontre de la m\^eme fa\c{c}on, en observant que pour tout $c>0$ fix\'e, on a
$$
n^{-cn}\zeta_{n,k}(-cn) \to \frac{e^c-2}{e^c-1}, \quad  (n>k \to \infty).
$$
C'est donc le point  \textit{(ii)} qui va retenir notre attention. Dans \cite{borwein}, Borwein, Fee, Ferguson et van der Waall d\'emontrent que $\varphi_p=\rho_{p,p-1}$ pour $p$ premier (Theorem 4.10). Leur approche peut \^etre \'etendue au cas de tous les entiers, en faisant appel au th\'eor\`eme d'\'equivalence de Bohr, que nous rappelons maintenant.

Deux s\'eries de Dirichlet ordinaires
$$ f(s) = \sum_{m=1}^\infty \frac{a(m)}{m^s}  \quad \mbox{ et } \quad 
g(s) = \sum_{m=1}^\infty \frac{b(m)}{m^s}$$
sont dites \'equivalentes s'il existe une fonction $f: \mathbb{N}^* \to \mathbb{C}$ compl\`etement multiplicative, telle que 
$b(m) = a(m) f(m)$ et $|f(m)|=1$ pour tout $m \soe 1$ (cf. \cite[Theorem 8.12]{apostol2}, o\`u la condition impos\'ee \`a $f$, apparemment plus faible, est en fait \'equivalente \`a la n\^otre). Le th\'eor\`eme d'\'equivalence de Bohr consiste alors en l'assertion suivante.
\begin{quote}
\textit{Soient $f(s)$ et $g(s)$ deux s\'eries de Dirichlet \'equi\-va\-lentes qui convergent absolument pour $\sigma > \sigma_0$. Alors, $f(s)$ et $g(s)$ prennent le m\^eme ensemble de valeurs dans toute bande ouverte incluse dans leur domaine de convergence absolue. Plus pr\'ecisement, si $\sigma_0 \ioe a < b$, alors
$$ \bigl\{ f(s), \,  a < \sigma < b \bigr\} = \bigl\{ g(s), \, a < \sigma < b \bigr\}.$$}
\end{quote}

Cet \'enonc\'e figure dans \cite[Theorem 8.16]{apostol2}, o\`u seul le cas des demi-plans $\sigma >a$ est signal\'e. La d\'emonstration pour une bande verticale est identique.

En particulier, $\zeta_n(s)$ prend dans tout demi-plan $\sigma <b$ exactement les m\^emes valeurs que 
$$
\zeta_{n,\chi}(s)=\sum_{m=1}^n \chi(m)m^{-s},
$$
o\`u $\chi$ est une fonction compl\`etement multiplicative \`a valeurs $\pm 1$. Si
\begin{equation} \label{chipm}
\chi(n) = \pm 1, \quad \chi(n-1) = \chi(n-2) = \dots = \chi(n-k) = \mp 1,
\end{equation}
on a 
$$
\pm \zeta_{n,\chi}(\sigma) \ioe \zeta_{n,k}(\sigma), \quad (\sigma \in \mathbb{R}).
$$
Par cons\'equent, $\zeta_{n,\chi}(\sigma)$ prend la valeur $0$ dans la demi-droite $]-\infty,\rho_{n,k}]$, et
donc $\zeta_n(s)$ poss\`ede un z\'ero dans le demi-plan $\sigma < b$ si $b>\rho_{n,k}$. 
Cela prouve \textit{(ii)}, qui repose donc sur la proposition suivante, dont la d\'emonstration est l'objet des paragraphes \ref{t2}, \ref{t3} et \ref{t7}.

\begin{prop} \label{kvalores} Soit $k \soe 1$. Il existe $n_1(k)$ dans $ \mathbb{N}^*$ tel que, pour tout $n \soe n_1(k)$, il existe 
une fonction $\chi: \mathbb{N}^* \to \{ \pm 1\}$ compl\`etement multiplicative v\'erifiant \eqref{chipm}.
\end{prop}

Nous verrons que l'outil essentiel de la d\'emonstration de la proposition \ref{kvalores} est le th\'eor\`eme de Siegel sur la finitude du nombre de points entiers d'une courbe elliptique d\'efinie sur les rationnels.

Terminons ce paragraphe en rappelant le r\'esultat de Montgomery \cite{montgomery} :
$$  \sup \bigl\{ \sigma \mid \zeta_n(s) = 0\}= 1 + \Bigl( \frac{4}{\pi} - 1 + o(1) \Bigr) \frac{\log \log n}{\log n}, \quad n \to +\infty, $$
dont la d\'emonstration utilise \'egalement le th\'eor\`eme d'\'equivalence de Bohr. 

\section{Syst\`eme d'\'equations diophantiennes quadratiques}\label{t2}

Nous utiliserons au \S \ref{t3} la propri\'et\'e de finitude suivante.
\begin{prop}\label{t4}
Soit $u_1>0$, $u_2>0$, $u_3>0$, $k_3>k_2>k_1 \soe 0$ des nombres entiers. Le nombre de triplets $(x,y,z)$ de nombres entiers tels que
\begin{equation}
  \label{t5}
  u_1x^2+k_1=u_2y^2+k_2=u_3z^2+k_3,
\end{equation}
est fini.
\end{prop}
\dem
  
Si $(x,y,z)$ est une solution de \eqref{t5}, alors $(X,Y)=(u_2u_3u_1^2x^2,u_1^2u_2^2u_3^2xyz)$ est solution de 
\begin{equation}
  \label{t6}
Y^2=X(X-\alpha)(X-\beta),  
\end{equation}
 avec $\alpha =u_1u_2u_3(k_2-k_1)$, $\beta = u_1 u_2 u_3 (k_3-k_1)$.
L'\'equation \eqref{t6} est celle d'une courbe elliptique sous la forme de Weierstrass. Comme $0<\alpha < \beta$, le th\'eor\`eme de Siegel \cite[Chapter IX, Corollary 3.2.1]{silverman} nous permet d'affirmer que \eqref{t6} n'a qu'un nombre fini de solutions enti\`eres $(X,Y)$. Par suite, \eqref{t5} n'a qu'un nombre fini de solutions enti\`eres $(x,y,z)$.\fin

\section{Partie sans facteur carr\'e d'entiers cons\'ecutifs}\label{t3}

Pour $n \soe 1$, notons $\rad(n)$ sa partie sans facteur carr\'e : c'est l'unique entier $b \soe 1$ sans facteur carr\'e tel que $n=a^2 b$, avec $a$ entier. La proposition suivante rassemble les propri\'et\'es de la fonction $\rad(n)$ utilis\'ees au \S \ref{t7}.
\begin{prop}\label{t1}
  Pour $n>k \soe 1$, on consid\`ere les nombres
  \begin{equation}
  \label{t0}
   \rad(n), \rad(n-1),\dots,\rad(n-k)
  \end{equation}

(i) Soit $p$ un nombre premier, $p>k$. Alors $p$ divise au plus un des nombres \eqref{t0}.

(ii) Si $n \soe k^2+k$, les nombres \eqref{t0} sont deux \`a deux distincts.

(iii) Il existe $n_0(k)$ tel que, pour $n \soe n_0(k)$, au plus deux d'entre les nombres \eqref{t0} ont tous leurs facteurs premiers $\ioe k$.
\end{prop}
\dem
  
\textit{(i)} Si $p|\rad (n-i)$ et $p|\rad (n-j)$ avec $0\ioe i<j\ioe k$, alors $p|(j-i)\ioe k$.

\textit{(ii)} Si $n-i=ua^2$, $n-j=ua'^2$ avec $u,a,a' \in \mathbb{N}^*$, $0\ioe i<j\ioe k$, alors
$$
k \soe j-i =u(a^2-a'^2)=u(a+a')(a-a')>ua \soe(ua^2)^{1/2}\soe(n-k)^{1/2},
$$
donc $n <k^2+k$.

\textit{(iii)}  Soient $0 \ioe k_1 < k_2 < k_3 \ioe k$ tels que les nombres $ u_i = \rad(n-k_i)$, $1 \ioe i \ioe 3$, 
poss\`edent tous leurs diviseurs premiers $ \ioe k$. Puisque $u_i$ 
est sans facteur carr\'e, $u_i \mid P$, o\`u $ P = \displaystyle \prod_{p \ioe k} p$ est le produit des premiers 
$ \ioe k$. On \'ecrit $ n-k_i = u_i a_i^2$, o\`u $a_i \in \mathbb{N}^*$, $1 \ioe i \ioe 3$. Pour chaque choix des $u_i$ et des $k_i$, l'ensemble des entiers $n$ correspondants est fini d'apr\`es la proposition \ref{t4}. Puisque l'ensemble des sextuples $(u_1,u_2,u_3,k_1,k_2,k_3)$ est fini (avec moins de $P^3 k^3$ \'el\'ements), cela d\'emontre le r\'esultat.\fin

\section{Fonctions compl\`etement multiplicatives \`a valeurs $\pm 1$ et entiers cons\'ecutifs}\label{t7}

On d\'emontre maintenant la proposition \ref{kvalores}. On commence par remarquer que $ \chi(m) = \chi \bigl( \rad(m) \bigr) $ pour tout $m \in  \mathbb{N}^*$ si $\chi$ est compl\`etement multiplicative \`a valeurs $\pm 1$.

Soit $n \soe n_1(k)$, o\`u $n_1(k)=\max \bigl (k^2+k,n_0(k)\bigr )$, et $n_0(k)$ v\'erifie le \textit{(iii)} de la proposition \ref{t1}.\\
Premier cas : on suppose qu'un nombre premier $p>k$ divise $\rad(n)$. D'apr\`es la proposition \ref{t1}, \textit{(i)}, les diviseurs premiers des $\rad(n-j)$, $0 < j \ioe k$, sont tous distincts de $p$. On pose alors $\chi(p)=-1$, et $\chi(q)=1$ pour $q \neq p$, pour obtenir 
(\ref{chipm}) avec $\chi(n)=-1$.\\
Deuxi\`eme cas : on suppose que tous les diviseurs premiers de $\rad(n)$ sont $\ioe k$. D'apr\`es la proposition \ref{t1}, \textit{(iii)}, il existe $j_0$ tel que $0<j_0 \ioe k$ et tel que $\rad(n-j)$ poss\`ede un facteur premier $p_j>k$ pour $0 < j \ioe k$, $j \not =j_0$. D'apr\`es la proposition \ref{t1}, \textit{(i)}, chacun des $p_j$ ne divise que $\rad(n-j)$ parmi les nombres \eqref{t0}. De plus, d'apr\`es la proposition \ref{t1}, \textit{(ii)}, on a $\rad(n) \neq \rad(n-j_0)$. Il existe donc un nombre premier $p$ qui ne divise qu'un seul de ces deux nombres. Premier sous-cas :  $p \mid \rad(n)$ et $p \nmid \rad(n-j_0)$. On pose alors $\chi(p)=-1$, $\chi(p_j)=(-1)^{[p \mid \rad(n-j)]}$ pour $0 < j \ioe k$, $j \not =j_0$, et $\chi(q)=1$ sur les autres nombres premiers. Ainsi (\ref{chipm}) est r\'ealis\'ee avec $\chi(n)=-1$. Deuxi\`eme sous-cas :  $p \nmid \rad(n)$ et $p \mid \rad(n-j_0)$. On pose alors $\chi(p)=-1$, $\chi(p_j)=(-1)^{[p \nmid \rad(n-j)]}$ pour $0 < j \ioe k$, $j \not =j_0$, et $\chi(q)=1$ sur les autres nombres premiers. Ainsi (\ref{chipm}) est r\'ealis\'ee avec $\chi(n)=1$.   

\smallskip

\noindent {\scriptsize {\bf Remerciements.}  Nous remercions Christian Ballot et Gary Walsh de nous avoir sugg\'er\'e l'utilisation du th\'eor\`eme de Siegel dans le traitement des syst\`emes d'\'equations diophantiennes.
Le premier auteur remercie le laboratoire Poncelet et l'Universit\'e Ind\'ependante de Moscou ; le second auteur remercie le LMNO et l'Universit\'e de Caen Basse-Normandie, 
et particuli\`erement Driss Essouabri, pour leur accueil et d'excellentes conditions de travail.}


\medskip

\begin{multicols}{2}
\footnotesize  

\noindent BALAZARD, Michel\\
Institut de Math\'ematiques de Luminy, UMR 6206\\
CNRS, Université de la Méditerranée\\
Case 907\\
13288 Marseille Cedex 09\\
FRANCE\\
Adresse \'electronique : \texttt{balazard@iml.univ-mrs.fr}

\smallskip

\noindent VEL\'ASQUEZ CASTA\~N\'ON, Oswaldo\\
Institut de Math\'ematiques de Bordeaux, UMR 5251\\
CNRS, Universit\'e Bordeaux 1\\
351, cours de la Lib\'eration\\
33405 Talence Cedex\\
FRANCE\\
Adresse \'electronique : \texttt{Oswaldo.Velasquez@math.u-bordeaux1.fr
}
\end{multicols}

\end{document}